\documentclass[12pt,a4paper,psamsfonts]{amsart}
\usepackage{amssymb,amscd,amsxtra,calc}
\usepackage{cmmib57}

\theoremstyle{plain}
    \newtheorem{thm}{Theorem}[section]

    \newtheorem{corollary}[thm]{Corollary}
    
    \newtheorem{proposition}[thm]{Proposition}
    \newtheorem{conjecture}[thm]{Conjecture}
    
    \newtheorem{theorem}[thm]{Theorem}

\theoremstyle{definition}

    \newtheorem{remark}[thm]{Remark}
\theoremstyle{remark}

    \newtheorem{setup}[thm]{}

\newcommand{\C}{\mathbb{C}}

\newcommand{\BPP}{\mathbb{P}}
\newcommand{\Q}{\mathbb{Q}}
\newcommand{\PP}{\mathbb{P}}

\newcommand{\R}{\mathbb{R}}
\newcommand{\BZZ}{\mathbb{Z}}
\newcommand{\Z}{\mathbb{Z}}

\newcommand{\OO}{\mathcal{O}}

\newcommand{\alb}{\operatorname{alb}}
\newcommand{\Alb}{\operatorname{Alb}}

\newcommand{\Aut}{\operatorname{Aut}}

\newcommand{\Exc}{\operatorname{Exc}}
\newcommand{\Gal}{\operatorname{Gal}}
\newcommand{\GL}{\operatorname{GL}}

\newcommand{\id}{\operatorname{id}}

\newcommand{\Ker}{\operatorname{Ker}}
\newcommand{\lin}{\operatorname{lin}}
\newcommand{\NE}{\overline{\operatorname{NE}}}

\newcommand{\Nklt}{\operatorname{Nklt}}

\newcommand{\NS}{\operatorname{NS}}

\newcommand{\rank}{\operatorname{rank}}

\newcommand{\Supp}{\operatorname{Supp}}
\newcommand{\torsion}{\operatorname{torsion}}

\begin{document}

\title[Birational geometry in the study of dynamics of automorphisms]
{Birational geometry in the study of dynamics of automorphisms and Brody/Mori/Lang hyperbolicity}

\author{De-Qi Zhang}
\address
{
\textsc{Department of Mathematics} \endgraf
\textsc{National University of Singapore} \endgraf
\textsc{10 Lower Kent Ridge Road, Singapore 119076
}}
\email{matzdq@nus.edu.sg}

\begin{abstract}
We survey our recent papers (some being joint ones) about the relation between the geometry of a
compact K\"ahler manifold and the existence of automorphisms of positive entropy on it. 
We also use the language of log
minimal model program (LMMP) in biraitonal geometry, but not its more sophisticated technical part.
We give applications of LMMP to positivity of log canonical divisor of a Mori / Brody / Lang
hyperbolic (quasi-) projective (singular) variety.
\end{abstract}



\thanks{This work is supported by an ARF of NUS
}

\maketitle

\tableofcontents

\section{Introduction}\

We work over the field $\C$ of complex numbers.
In \S 2, we consider automorphisms or automorphism groups of positive (topological)
entropy on compact K\"ahler manifolds or normal projective varieties.
This part is the continuation of our earlier survey paper \cite{ICCM07}.
In \S 3, we consider $\Aut(X)$ for an almost homogeneous compact K\"ahler manifold $X$.
In \S 4, we consider positivity of log canonical divisor
of Mori hyperbolic or Brody hyperbolic quasi-projective varieties.
In \S 5, we consider positivity of canonical divisor of algebraic Lang hyperbolic singular varieties.

We use the convention in \cite{DS}, \cite{KM}, and the book of Hartshorne.

\section{Automorphisms of positive entropy}

This section is mainly based on \cite{NullG}, \cite{maxrank} and the joint paper \cite{CWZ}.
We refer to Dinh \cite{Di} for another survey on results of Tits alternative type.

For a linear transformation $L$ on a finite-dimensional vector space $V$
over $\C$ or its subfields,
its {\it spectral radius} is defined as
$$\rho(L) := \max \{|\lambda| \, : \, \lambda \, \in \, \C \,\,\, \text{is an eigenvalue of} \,\,\, L\}.$$

Let $X$ be a compact complex K\"ahler manifold and $Y$ a normal projective variety,
and let $g \in \Aut(X)$ and $f \in \Aut(Y)$.
Define the (topological) {\it entropy} $h(*)$ and {\it first dynamical degrees} $d_1(*)$ as:

$$\begin{aligned}
h(g) : &= \log \rho(g^* \, | \, \oplus_{i \ge 0} \, H^i(X, \C)), \\
d_1(g) : &= \rho(g^* \, | \, H^2(X, \C)) \,\, (= \rho(g^* \, | \, H^{1,1}(X))), \\
d_1(f) : &= \rho(f^* \, | \, \NS_{\C}(Y))
\end{aligned}$$
where $\NS_{\C}(Y) := \NS(Y) \otimes_{\Z} \C$ is the {\it complexified Neron-Severi group}.
By the fundamental work of Gromov and Yomdin, the above definition of entropy is equivalent to its original definition
(cf. \cite[\S 2.2]{DS} and the references therein).
Further, when $Y$ is smooth, the above two definitions of $d_1(*)$ coincide;
for $\Q$-factorial $Y$ (cf. \cite[0.4(1)]{KM}), we have $d_1(f) = d_1(\widetilde{f})$
where $\widetilde{f}$ is the lifting of $f$ to the one on an $\Aut(Y)$-equivariant resolution of $Y$.
We call $\tau := g$ or $f$, of {\it positive entropy} (resp. {\it null entropy}) if
$d_1(\tau) > 1$ (resp. $d_1(\tau) = 1$), or equivalently
$h(\tau) > 0$ (resp. $h(\tau) = 0$) in the case of compact K\"ahler manifold.

We say that the induced action $G \, | \, H^{1,1}(X)$ is
{\it Z-connected} if its Zariski-closure in $\GL(H^{1,1}(X))$ is connected
with respect to the Zariski topology;
in this case, the {\it null set}
$$N(G) := \{g \in G \, | \, g \,\, \text{is of null entropy} \}$$
is a (necessarily normal) subgroup of $G$ (cf. \cite[Theorem 1.2]{Z-Tits}).
In \cite{Z-Tits}, we have proved:

\begin{theorem} (cf.~\cite{Z-Tits}) \label{Z-TitsTh}
Let $X$ be an $n$-dimensional $(n \ge 2)$
compact complex K\"ahler manifold and $G$ a subgroup of $\Aut(X)$.
Then one of the following two assertions holds:
\begin{itemize}
\item[(1)]
$G$ contains a subgroup isomorphic to the non-abelian free group $\BZZ * \BZZ$,
and hence $G$ contains subgroups isomorphic to non-abelian free groups of all countable ranks.
\item[(2)]
There is a finite-index subgroup $G_1$ of $G$ such that
the induced action $G_1 \, | \, H^{1,1}(X)$ is solvable and Z-connected.
Further, the subset
$$N(G_1) := \{g \in G_1 \,\, | \,\, g \,\,\, \text{\rm is of null entropy}\}$$
of $G_1$ is a normal subgroup of $G_1$ and the quotient group $G_1/N(G_1)$ is a
free abelian group of rank $r \le n-1$. We call this $r$ {\bf the rank of $G_1$ and denote it as $r = r(G_1)$}.
\end{itemize}
\end{theorem}

A compact K\"ahler manifold $X$ is {\it ruled}
if it is bimeromorphic to a manifold with
a $\BPP^1$-fibration.
By a result of Matsumura, $X$ is ruled if $\Aut_0(X)$ is not a compact torus
(cf. \cite[Proposition 5.10]{Fu}).
When $X$ is a compact complex K\"ahler manifold (or a normal projective variety),
set $L := H^2(X, \Z)/(\torsion)$ $($resp.\ $L := \NS(X)/(\torsion))$,
$L_{\R} := L \otimes_{\Z} \R$, and $L_{\C} := L \otimes_{\Z} \C$.

In \cite{CWZ}, we have strengthened Theorem \ref{Z-TitsTh} to the following:

\begin{theorem} (cf.~\cite[Theorem 1.5]{CWZ}) \label{ThA11}
Let $X$ be a compact K\"ahler $($resp.\ projective$)$ manifold
of dimension $n$ and $G \le \Aut(X)$ a subgroup.
Then one of the following properties holds.
\begin{itemize}
\item[(1)]
$G | L_{\C} \ge \Z * \Z$ $($the non-abelian free group of rank two$)$,
and hence $G \ge \Z * \Z$.
\item[(2)]
$G | L_{\C}$ is virtually solvable and
$$G \ge K \cap L(\Aut_0(X)) \ge \Z * \Z$$
where
$L(\Aut_0(X))$ is the linear part of
$\Aut_0(X)$ {\rm (cf.~\cite[Def.~3.1, p. 240]{Fu})} and
$$K = \Ker(G \to \GL(L_{\C})) ,$$
so $X$ is ruled
{\rm (cf. \cite[Proposition 5.10]{Fu})}.
\item[(3)]
There is a finite-index solvable subgroup $G_1$ of $G$
such that the null subset $N(G_1)$ of $G_1$
is a normal subgroup of $G_1$ and
$$G_1/N(G_1) \cong \Z^{\oplus r}$$
for some $r \le n - 1$.
\end{itemize}
In particular, either $G \ge \Z * \Z$ or $G$ is virtually solvable.
In Cases $(2)$ and $(3)$ above, $G | L_{\C}$ is finitely generated.
\end{theorem}

In view of Theorem \ref{Z-TitsTh}, we are more interested in the group $G \le \Aut(X)$ where
$G \, | \, H^{1,1}(X)$ is solvable and $Z$-connected and
that the rank $r(G) = \dim X - 1$ (maximal value).

In the following, denote by $\Aut_0(X)$ the {\it identity connected component} of $\Aut(X)$.
A group {\it virtually has a property} (P) if a finite-index subgroup of it has the property (P).
Let $X$ be a normal projective variety with at worst canonical singularities (cf. \cite[Definition 2.34]{KM})
and $G \le \Aut(X)$ a subgroup.

The pair $(X, G)$ is {\it non-minimal} if:
there are a finite-index subgroup $G_1$ of $G$ and a {\it non-isomorphic} $G_1$-equivariant birational morphism $X \to Y$
onto a normal projective variety $Y$ with at worst {\it isolated} canonical singularities. The pair
$(X, G)$ is {\it minimal} if it is not non-minimal.

A complex torus has lots of symmetries.
Conversely, our Theorem \ref{ThA1} in \cite{NullG} (see also Theorem \ref{ThC1} for non-algebraic manifolds)
says that the maximality $r(G) = \dim X - 1$
occurs only when $X$ is a quotient of a complex torus $T$
and $G$ is mostly descended from the symmetries on the torus $T$.

%

\begin{theorem} (cf. \cite[Theorem 1.1]{NullG})\label{ThA1}
Let $X$ be an $n$-dimensional $(n \ge 3)$
normal projective variety and $G \le \Aut(X)$ a subgroup such that
the induced action $G \, | \, \NS_{\C}(X)$ is solvable and $Z$-connected and that
the rank $r(G) = n - 1$ $($i.e., $G/N(G) = \Z^{\oplus n-1})$.
Assume the
three conditions:

\begin{itemize}
\item[(i)]\footnote{In Theorem 1.2, the condition (i) can be weakened to:
(i)' X has at worst canonical singularities and is smooth in codimension two.
This is by the same proof of ours and due to Theorem 1.16 in the recent paper
of D. Greb, S. Kebekus and T. Peternell:
`\'Etale fundamental groups of Kawamata log terminal spaces, flat sheaves, and quotients of Abelian varieties', arXiv:1307.5718.
}
$X$ has at worst canonical, quotient
singularities {\rm (cf. \cite[Definition 2.34]{KM})}.
\item[(ii)]
$X$ is a minimal variety, i.e., the canonical divisor $K_X$ is nef {\rm (cf. \cite[0.4(3)]{KM})}.
\item[(iii)]
The pair $(X, G)$ is minimal in the sense above.
\end{itemize}
Then the following four assertions hold.
\begin{itemize}
\item[(1)]
The induced action $N(G) \, | \, \NS_{\C}(X)$ is a finite group.
\item[(2)]
$G \, | \, \NS_{\C}(X)$ is a virtually free abelian group of rank $n-1$.
\item[(3)]
Either $N(G)$ is a finite subgroup of $G$ and hence
$G$ is a virtually free abelian group of rank $n-1$, or
$X$ is an abelian variety and the group
$N(G) \cap \Aut_0(X)$ has finite-index in $N(G)$ and is Zariski-dense in $\Aut_0(X)$ $(\cong X)$.
\item[(4)]
We have $X \cong T/F$ for a finite group $F$ acting freely outside a finite set of an abelian variety $T$.
Further,
for some finite-index subgroup $G_1$ of $G$, the action of $G_1$ on $X$ lifts to
an action of $G_1$ on $T$.
\end{itemize}
\end{theorem}

When the $X$ below is non-algebraic, we don't require the minimality of the pair $(X, G)$.

\begin{theorem} (cf.~\cite[Theorem 2.1]{maxrank}) \label{ThC1}
Let $X$ be an $n$-dimensional $(n \ge 3)$
compact complex K\"ahler manifold which is not algebraic.
Let $G \le \Aut(X)$ be a subgroup such that
the induced action $G \, | \, H^{1,1}(X) \le \Aut(H^{1,1}(X))$ is solvable and $Z$-connected and
that the rank $r(G) = n - 1$ $($i.e., $G/N(G) = \Z^{\oplus n-1})$.
Assume that $X$ is minimal, i.e., the canonical divisor $K_X$ is
contained in the closure of the K\"ahler cone of $X$.
Then the following four assertions hold.
\begin{itemize}
\item[(1)]
The induced action $N(G) \, | \, H^{1,1}(X)$ is a finite group.
\item[(2)]
$G \, | \, H^{1,1}(X)$ is a virtually free abelian group of rank $n-1$.
\item[(3)]
Either $N(G)$ is a finite subgroup of $G$ and hence
$G$ is a virtually free abelian group of rank $n-1$, or
$X$ is a complex torus and the group
$N(G) \cap \Aut_0(X)$ has finite-index in $N(G)$ and is Zariski-dense in $\Aut_0(X)$ $(\cong X)$.
\item[(4)]
We have $X \cong T/F$ for a finite group $F$ acting freely outside a finite set of a complex torus $T$.
Further,
for some finite-index subgroup $G_1$ of $G$, the action of $G_1$ on $X$ lifts to
an action of $G_1$ on $T$.
\end{itemize}
\end{theorem}

In Theorem \ref{ThA2} below, we will assume that (i) $G$ is abelian and (ii) the absence
of point-wise $G$-fixed subvarieties of positive dimension or $G$-periodic rational curves or $Q$-tori.
In Theorem \ref{ThA1}, these two restrictions are replaced by the natural minimality condition on $X$
and the pair $(X, G)$, and that $G \, | \, \NS_{\C}(X)$ is solvable, the latter of which is natural
in view of Theorem \ref{Z-TitsTh}. The quotient singularities assumption in Theorem \ref{ThA1}
is necessary because an effective characterization of torus quotient is only available in dimension three
by \cite{SW} where the bulk of the argument is to show that the variety has only quotient singularities.

The lack of the abelian-ness assumption on $G$ makes our argument much harder, for instance we cannot
simultaneously diagonalize $G \, | \, \NS_{\C}(X)$ or find enough number of linearly independent
common nef eigenvectors of $G$ as required in \cite{DS} for abelian groups.

Theorems \ref{ThA1} and \ref{ThC1} answer \cite[Question 2.17]{Z-Tits}, assuming the conditions here.
When $G$ is abelian, the finiteness of $N(G)$ is proved in the inspiring paper
of Dinh-Sibony \cite[Theorem 1]{DS} (cf. also \cite{CY3}),
assuming only $r(G) = n-1$.
For non-abelian $G$, the finiteness of $N(G)$ is not true and
we can at best expect
that $N(G)$ is virtually included in $\Aut_0(X)$ (as done in Theorems \ref{ThA1} and \ref{ThC1}),
since a larger group $\widetilde{G} := \Aut_0(X) \, G$ satisfies
$$\widetilde{G} \, | \, \NS_{\C}(X) = G \, | \, \NS_{\C}(X), \,\,
N(\widetilde{G}) = \Aut_0(X) . N(G) \ge \Aut_0(X), \,\,
\widetilde{G}/N(\widetilde{G}) \cong G/N(G).$$

There are examples $(X, G)$ with rank $r(G) = \dim X - 1$ and $X$ complex tori or their quotients
(cf. \cite[Example 4.5]{DS}, \cite[Example 1.7]{CY3}).

The proof of Theorem \ref{ThA1} is much harder than that of Theorem \ref{ThC1}
because of the presence of singularities on $X$.

The conditions (i) - (iii) in Theorem \ref{ThA1} are quite necessary in deducing
$X \cong T/F$ as in Theorem \ref{ThA1}(4).
Indeed, if $X \cong T/F$ as in Theorem \ref{ThA1}(4), then $X$ has only quotient singularities
and $d K_X \sim 0$ (linear equivalence) with $d = |F|$, and we may even assume that $X$ has only canonical singularities
if we replace $X$ by its global index-$1$ cover; thus $X$ is a minimal variety.
If the pair $(X, G)$ is not minimal so that there
is a non-isomorphic $G_1$-equivariant birational morphism
$X \to Y$ for some finite-index subgroup $G_1$ of $G$, then the exceptional locus of this morphism
is $G_1$- and hence $G$-periodic, contradicting the fact that the rank $r(G) = n-1$.

The {\it first key step} in proving Theorem \ref{ThA1} is the
analysis of our quasi-nef sequence $L_1 \dots L_k$ ($0 < k < n$)
(cf. \cite[\S 2.2]{Z-Tits}) and we are able to show that $L_i$ can actually
be taken to be nef, when the rank $r(G) = n-1$.
The {\it second key step} is
to split $G$ as $N(G) \, H$
such that $H \, | \, H^{1,1}(X)$ is free abelian; thus we have a nef and big class $A$
as the sum of nef common eigenvectors of $H$,
leading to the vanishing of $A^{n-i} . c_i(X)$,
where $c_i(X)$ ($i = 1, 2$) are Chern classes (cf. \cite[pages 265-267]{SW}).
The {\it third key step} is to use the minimality of $(X, G)$ and Kawamata's base point freeness
result for $\R$-divisor (cf. \cite[Theorem 3.9.1]{BCHM}, available only for projective
variety at the moment) to deduce the vanishing
of $c_2(X)$ as a linear form; this vanishing does not directly follow from the vanishing of $A^{n-2} . c_2(X)$
because $A$ may not be ample. Now Theorem \ref{ThA1}(4) follows from the vanishing of $c_i(X)$ ($i = 1, 2$)
and the characterization of torus quotient originally deduced from Yau's deep result
(cf. \cite{Be}).

\begin{remark}\label{rThA1}
$ $
\begin{itemize}
\item[(1)]
When $\dim X = 3$, the (i) in Theorem \ref{ThA1} can be replaced by:
(i)' $X$ has at worst canonical singularities.
See \cite[Corollary at p. 266]{SW}.
\item[(2)]
Theorems \ref{ThA1} and \ref{ThC1} are not true when $n := \dim X = 2$;
see Remark \ref{rThA2} below.
We used $n \ge 3$ to deduce the vanishing of $c_2(X) . A^{n-2}$
as commented above.
\end{itemize}
\end{remark}

For the definitions of {\it Kodaira dimension}
$\kappa(X)$ and singularities of {\it terminal}, {\it canonical} or {\it klt type},
we refer to \cite[Definitions 2.34 and 7.73]{KM}.
A subvariety $Y \subset X$ is called $G$-{\it periodic}
if $Y$ is stabilized (set theoretically) by a finite-index subgroup of $G$.
Denote by $q(X) := h^1(X, \OO_{X})$ the {\it irregularity} of $X$.
A projective manifold $Y$ is a $Q$-{\it torus}
if $Y = A/F$ for a finite group $F$ acting freely on an abelian variety $A$.

Without assuming the minimality of the pair $(X, G)$
we have proved in \cite{maxrank}:

\begin{theorem} (cf.~\cite[Theorem 1.2]{maxrank}) \label{ThA2}
Let $X$ be an $n$-dimensional normal projective variety
and let $G := {\Z}^{\oplus n-1}$ act on $X$ faithfully
such that every element of $G \setminus \{\id\}$ is of positive entropy.
Then the following hold:

\begin{itemize}
\item[(1)]
Suppose that
$\tau : A \to X$ is a $G$-equivariant
finite surjective morphism from an abelian variety $A$. Then $\tau$
is \'etale outside a finite set $($hence $X$ has only quotient
singularities and is $klt)$; $K_X \sim_{\Q} 0$ $(\Q$-linear equivalence$)$;
no positive-dimensional proper subvariety $Y \subset X$ is $G$-periodic.

\item[(2)]
Conversely, suppose that
no positive-dimensional $G$-periodic
subvariety $Y \subset X$ is either fixed $($point wise$)$ by a finite-index subgroup of $G$,
or is a $Q$-torus with $q(Y) > 0$,
or has $\kappa(Y) = -\infty$.
Suppose
also one of the following two conditions.

\item[(2a)]
$n = 3$, and $X$ is $klt$.
\item[(2b)]
$n \ge 3$, and $X$ has only quotient singularities.

Then $X \cong A/F$ for a finite group $F$ acting freely outside a finite set of an abelian variety $A$.
Further, for some finite-index subgroup $G_1$ of $G$, the action of $G_1$ on $X$ lifts to
an action of $G_1$ on $A$.
\end{itemize}
\end{theorem}

As a consequence of Theorems \ref{ThA2} and \cite[Theorem 1]{DS} (or Theorem \ref{Z-TitsTh}), we have:

\begin{corollary} (cf.~\cite{maxrank}) \label{CorB}
Let $X$ be a normal projective variety of dimension $n \ge 3$ with
only
quotient singularities,
and let $G := {\Z}^{\oplus r}$ act on $X$ faithfully for some
$r \ge n-1$
such that every element of $G \setminus \{\id\}$ is of positive entropy.
Then $r = n-1$. Further,
$X \cong A/F$ for a finite group $F$ acting freely
outside a finite set of an abelian variety $A$,
if and only if
$X$ has no positive-dimensional $G$-periodic proper subvariety $Y \subset X$.
\end{corollary}

The proof of \cite[Theorem 2.2]{maxrank} gives the following, essentially already proved
in \cite{CY3}.

\begin{corollary} (cf.~\cite{maxrank}) \label{CorD}
Let $X$ be a normal projective variety of dimension $n \ge 3$
and let $G := {\Z}^{\oplus r}$ act on $X$ faithfully for some
$r \ge n-1$
such that every element of $G \setminus \{\id\}$ is of positive entropy.
Then $r = n-1$.
Suppose further that both $X$ and $(X, G)$ are minimal in the sense of \cite[Definition 2.1]{maxrank}
(which is slightly different from the one given before Theorem \ref{ThA1})
and either $X$ has only quotient singularities, or $X$ is a $klt$ threefold.
Then $X \cong A/F$ for a finite group $F$ acting freely
outside a finite set of an abelian variety $A$.
\end{corollary}

\begin{remark}\label{rThA2}
(1) In Theorem \ref{ThA2} (2)
we need to assume that $\dim X \ge 3$
which is used at the last step of the proof to show the vanishing of
the second Chern class $c_2(X)$.
In fact, inspired by the comment of the referee of \cite{maxrank}, one notices that
a complete intersection $X \subset \PP^2 \times \PP^2$
of two very general hypersurfaces of type $(1, 1)$ and $(2, 2)$
is a $K3$ surface (called Wehler's surface) of Picard number two,
$X$ has an automorphism $g$ of entropy $2 \log(2 + \sqrt{3}) > 0$, and
$X$ contains no $(-2)$-curve, so there is no $g$-periodic curve on $X$.
Thus, both $X$ and the pair $(X, \langle g \rangle)$ are minimal in the sense of
\cite[Definition 2.1]{maxrank}.
However, $X$ is not birational to the quotient of a complex torus,
because a (smooth) projective $K3$ surface $X$ birational to the quotient of a complex torus
has the transcendental lattice of rank $\le$
(that of a complex $2$-torus), i.e., $\le 5$, and hence has Picard number $\ge (h^2(X, \C) - 5)$
which is $17$.

(2) In Theorem \ref{ThA2} (2)
we can weaken the assumption on $G$
as a condition on $G^* := G_{|\NS_{\R}(X)}$ (cf. \cite[Theorem 4.7]{DS}):

\par \vskip 1pc \noindent
 ``$G^* \cong \Z^{\oplus n-1}$
and every element of $G^* \setminus \{\id\}$
is of positive entropy."

\par \vskip 1pc \noindent
But we need also to replace the last sentence

\par \vskip 1pc \noindent
``Further, \dots of $G_1$ on $A$." in Theorems \ref{ThA2} (2)
as:

\par \vskip 1pc \noindent
``Further, the action of $G$ on $X$ lifts to an action of a group
$\widetilde{G}$ on $A$ with $\widetilde{G}/\Gal(A/X) \cong G$."

\end{remark}

Our bimeromorphic point of view, in terms of the minimality assumption in
\ref{CorD}
towards the dynamics study seems natural,
since one may blow up some Zariski-closed and $G$-stable proper subset of $X$
(if such subset exists) to get another pair
$(X', G)$ which is essentially the same as the original pair $(X, G)$.

The very starting point of our proof is the existence of enough
nef eigenvectors $L_i$ of $G$,
due to the fundamental work of Dinh-Sibony \cite{DS}. The minimal model program
(cf. \cite{KM}) is used with references provided for non-experts. Our main contribution lies
in
showing that the pair
$(X, G)$ can be replaced with
an equivariant one so that $H:= \sum_{i=1}^n \, L_i$ is an ample divisor.
To conclude, we prove a result of Hodge-Riemann type for singular varieties
to show the vanishing of Chern classes $c_i(X)$ ($i = 1, 2$),
utilizing $H^{n-i} \cdot c_i(X) = 0 = H^{n-2} \cdot c_1(X)^2$.
Then we use the characterization of \'etale quotient of a complex torus as
the compact K\"ahler manifold $X$ with vanishing Chern classes $c_i(X)$ ($i = 1, 2$)
(cf. \cite[\S 1]{Be}) and its generalization to singular varieties (cf. \cite{SW}).

For a possible generalization of the proof to a K\"ahler $n$-fold $X$,
we remark that the restriction `$\rank r(G) = n-1$' implies that $X$ is either Moishezon and hence projective,
or has algebraic dimension $a(X) = 0$ (cf. \cite[Theorem 1.2]{Z-Tits}).
Thus the case $a(X) = 0$ remains to be treated. See a related remark in
\cite[\S 3.6]{Ca06}.

\section{Automorphism groups of almost homogeneous varieties}

This section is based on the joint paper \cite{FZ}.
Let $X$ be a compact K\"ahler manifold. Denote by $\Aut(X)$ the automorphism group of $X$ and by $\Aut_0(X)$ the identity connected component of $\Aut(X)$. By \cite{Fu}, $\Aut_0(X)$ has a natural meromorphic group structure.
Further there exists a unique meromorphic subgroup, say $L(X)$, of
$\Aut_0(X)$, which is meromorphically isomorphic to a linear algebraic group and such
that the quotient $\Aut_0(X)/L(X)$ is a complex torus.
In the following, by a subgroup of $\Aut_0(X)$ we always mean a meromorphic subgroup and by a linear algebraic subgroup
of $\Aut_0(X)$ we mean a Zariski closed meromorphic subgroup contained in $L(X)$.

For a subgroup $G \le \Aut(X)$, the
pair $(X, G)$ is called {\em strongly primitive} if for
every finite-index subgroup $G_1$ of $G$,
$X$ is not bimeromorphic to a non-trivial $G_1$-equivariant fibration,
i.e., there does not exist any compact K\"ahler manifold
$X'$ bimeromorphic to $X$, such that $X'$ admits  a
$G_1$-equivariant holomorphic map $X' \to Y$ with $0 < \dim Y < \dim X$  and  $G_1 \le \Aut(X')$.
See Oguiso - Truong \cite{OT} for the very first examples of strongly primitive pairs $(X, \langle g \rangle)$
with $X$ a Calabi-Yau threefold or rational threefold.
From the dynamical point of view, these manifolds are essential.
Our Theorem \ref{ThA3} of \cite{FZ} says that for these manifolds, unless it is  a complex torus,
there is no interesting dynamics if its automorphism group has positive dimension.

\begin{theorem} (cf.~\cite[Theorem 1.1]{FZ}) \label{ThA3}
Let $X$ be a compact K\"ahler manifold and $G \le \Aut(X)$ a
subgroup of automorphisms. Assume the following three conditions.
\begin{itemize}
\item[(1)] $G_0 := G \cap \Aut_0(X)$ is infinite.
\item[(2)]  $|G : G_0| = \infty $.
\item[(3)] The pair $(X, G)$ is strongly primitive.
\end{itemize}
Then $X$ is a complex torus.
\end{theorem}

As a key step towards Theorem \ref{ThA3}, we prove the following result in \cite{FZ}.
A proof for Theorem \ref{ThB3}(2) is long overdue
(and we do it geometrically via \ref{ThB3}(1)), but the authors could not find
it in any literature, even after consulting many experts across the continents.

\begin{theorem} (cf.~\cite[Theorem 1.2]{FZ}) \label{ThB3}
Let $X$ be a compact K\"ahler manifold and $G_0 \subset \Aut_0(X)$ a linear algebraic subgroup.
Assume that $G_0$ acts on $X$ with a Zariski open dense orbit.
Then we have:
\begin{itemize}
\item[(1)]
$X$ is projective; the anti canonical divisor $-K_X$ is big, i.e.
$\kappa(X, -K_X) = \dim X$.
\item[(2)]
$\Aut(X)/\Aut_0(X)$ is finite.
\end{itemize}
\end{theorem}

\begin{remark}\label{rThm3}
(i) The condition (2) in Theorem \ref{ThA3} is satisfied if  $G$ acts on  $H^2(X, \C)$
as an infinite group, or if $G$ has an element of positive entropy.

(ii)  Theorem \ref{ThB3} implies that when $\dim X \ge 3$
the case(4) in \cite[Theorem 1.2]{Z-Tits} does not occur, hence it can be removed from the statement.

(iii) Theorem \ref{ThA3} generalizes \cite[Theorem 1.2]{CWZ}, where it is proven for $\dim X = 3$ and under an additional assumption.

(iv) A similar problem
for endomorphisms of
homogeneous varieties has been studied by S. Cantat in \cite{Ca04}.
\end{remark}


Two applications are given. The first one
generalizes  the following result due to Harbourne  \cite[Corollary (1.4)]{Hb} to higher dimension:
if $X$ is a smooth projective rational surface with $\Aut_0(X) \ne (1),$
then $\Aut(X)/\Aut_0(X)$ is finite. To state it, recall (\cite[Theorem 4.1]{Fu})
that for any connected subgroup $H \le \Aut(X)$, there
exist a quotient space $X/H$ and an $H$-equivariant dominant meromoprhic map $X \dasharrow X/H$,
which satisfies certain universal property.

\begin{corollary} (cf.~\cite[1.4]{FZ})
Let $X$ be a compact K\"ahler manifold with irregularity $q(X) = 0$.
Suppose that the quotient space $X/\Aut_0(X)$
has dimension $\le 1$. Then $X$ is projective and $\Aut(X)/\Aut_0(X)$ is finite.
\end{corollary}

The second application essentially says that when we study dynamics of a compact K\"ahler manifold $X$,
we may assume that $\Aut_0(X)_{\lin} = (1)$, where $\Aut_0(X)_{\lin}$ is the largest connected linear
algebraic subgroup of $\Aut_0(X)$.

\begin{corollary} (cf.~\cite[1.5]{FZ}) \label{ThC3}
Let $X$ be a smooth projective variety and $G_0 \lhd G \le \Aut(X)$.
Suppose that $G_0$ is a connected linear closed subgroup of $\Aut_0(X)$.
Let $Y$ be a $G$-equivariant resolution of the quotient space
$X/G_0$ and replace $X$ by a $G$-equivariant resolution
so that the natural map $\pi : X \to Y$
is holomorphic.
Then for any $g \in G$,
we have the equality of the first dynamical degrees:
$$d_1(g_{| X}) = d_1(g_{| Y}),$$
where $d_1(g_{| X}) := \max\{|\lambda| \, ; \, \lambda \,\, \text{is
an eigenvalue of} \,\, g^* \, | \, H^{1,1}(X) \}$.

In particular, $G_{| X}$ is of null entropy if and only if so is $G_{| Y}$.
\end{corollary}

\begin{remark}
If $q(X)=0$, then $\Aut_0(X)$ (hence $G_0$) is always a linear algebraic group. On the other hand,
if $G_0$ is not linear, then Corollary \ref{ThC3} does not hold.
For example, let $X = T_1 \times T_2$ be the product of two complex tori and $G_0 = \Aut_0(T_1)$. The quotient space $Y:=X/G_0$ is $T_2$. Suppose that $g \in \Aut(T_1)$ is an element of positive entropy
 which acts trivially on $T_2$. Then we have
 $d_1(g_{| X}) > d_1(g_{| Y}) = 1$.
\end{remark}

\section{Positivity of log canonical divisors and Mori/Brody hyperbolicity}

This section is based on the results in the joint paper \cite{LZ}.
An algebraic variety $X$ is called {\it Brody hyperbolic} = BH
(resp.~Mori hyperbolic = MH)
if the following hypothesis (BH) (resp.~(MH)) is satisfied:

\begin{itemize}
\item [(BH)]
Every holomorphic map from the complex line $\C$ to $X$ is a constant map.
\item[(MH)]
Every algebraic morphism from the complex line $\C$ to $X$ is a constant map.\\
(Eqivalently, no algebraic curve in $X$ has normalization equal to
$\PP^1$ or $\C$.)
\end{itemize}

Clearly, (BH) $\Rightarrow$ (MH).

Consider a pair $(X, D)$ of an (irreducible) algebraic variety $X$ and a
(not necessarily irreducible
or equi-dimensional) Zariski-closed subset $D$ of $X$.
When $X$ is an irreducible curve and $D$ a finite subset of $X$,
the pair $(X, D)$ is {\it Brody hyperbolic}
(resp.~MH)
if $X \setminus D$ is
Brody hyperbolic (resp.~MH).
Inductively, for an algebraic variety $X$
and a Zariski-closed subset $D$ of $X$,
with $D = \cup_i D_i$ the irreducible decomposition, the pair $(X, D)$ is
{\it Brody hyperbolic} (resp.~MH)
if
$X \setminus D$ and the pairs $(D_k, D_k \cap (\cup_{i \ne k} D_i))$ for all $k$, are so.
Note that the Zariski-closed subset $D_k \cap (\cup_{i \ne k} D_i)$ may not be equi-dimensional in general.

Here, it is crucial from the perspective of hyperbolic geometry
that $D$ is a divisor and all the $D_k$'s are
at least $\Q$-Cartier
(so as to have equidimensionality above
for example). But we will work with
a bit less from the outset as demanded by our general setup.

The pair $(X, D)$ is called {\it projective} if $X$ is a projective variety and {\it log smooth} if
further $X$ is smooth and $D$ is a reduced {\it divisor} with simple normal crossings.\\[-3mm]

{\it From now on,
we always assume that $D$ is a 
reduced Weil divisor for a pair $(X, D)$.} \\[-3mm]

Recall that a divisor $F$ on a smooth projective variety
$X$ is called {\it numerically effective}
(nef), if $F.C\geq 0$ for all curves $C$ on $X$.
We consider the following conjecture:

\begin{conjecture}\label{ConjA}
Let $(X, D)$ be a log smooth projective Mori hyperbolic
$($resp.~Brody hyperbolic$)$ pair.
Then the log canonical divisor $K_X + D$ is nef
$($resp.~ample$)$.
\end{conjecture}
%

Our main theorems of \cite{LZ} below give an
affirmative answer to Conjecture \ref{ConjA} but with
further assumptions for the ampleness of $K_X + D$,
such as the ampleness of
at least $n-2$ irreducible
Cartier components of $D$ for
$n=\dim\,(X, D) := \dim X$. To precise them,
we will need to refer to the following two
standard conjectures on the structure of algebraic varieties.

\begin{conjecture}\label{ConjAA}{\rm\bf Abundance($l$):}
Let $(X, D)$ be an $l$-dimensional log smooth projective pair
whose log canonical divisor $K_X + D$ is nef.
Then some positive multiple of $K_X+D$ is base point free,
i.e. $K_X+D$ is semi-ample.
\end{conjecture}

\begin{conjecture}\label{ConjAAA} {\rm\bf CY($m$):}
Let $X$ be an $m$-dimensional simply connected
nonsingular projective variety with trivial canonical line bundle,
i.e. a Calabi-Yau manifold.
Then $X$ contains a rational curve.
\end{conjecture}

We remark that {\rm  Abundance($l$)} is known to hold
for $l \leq 3$ (even for dlt or lc pairs)
and {\rm CY($m$)} is known to hold
for $m\leq 2$ (see \cite[\S 3.13]{KM}, \cite{MM}).
\\[-3mm]

Let $D$ be a reduced divisor on a projective variety $X$
and $\mathcal D$ the collection of terms in its irreducible
decomposition $D = \sum_{i=1}^s D_i$. A stratum
of $D$ is a set of the form
$\cap_{i \in I} D_i$ ($= X$ when $I = \emptyset$)
for some partition $\{1, \dots, s\} = I \coprod J$.
We call a $\mathcal D$-{\it rational curve} a
rational curve $\ell$ in a stratum $\cap_{i \in I} D_i$ of $D$
such that the normalization of $\ell \setminus (\cup_{j \in J} D_j)$
contains the complex line $\C$. We call a $\mathcal D$-{\it algebraic $1$-torus}
a rational curve $\ell$ in a stratum $\cap_{i \in I} D_i$ of $D$ such that
$\ell \setminus (\cup_{j \in J} D_j)$ has the
$1$-dimensional algebraic torus
$\C^* = \C \setminus \{0\}$
as its normalization.
We call a closed subvariety $W\subseteq X$ a
$\mathcal{D}$-{\it compact variety}
if $W$ lies in some stratum $\cap_{i \in I} D_i$ of $D$
and $W\cap \cup_{j \in J} D_j=\emptyset$.
A $\mathcal{D}$-compact variety $W \subseteq X$ is called a
$\mathcal{D}$-compact rational curve,
a $\mathcal{D}$-{\it torus},
a $\mathcal{D}$-$\Q$-{\it torus}, or a $\mathcal{D}$-$\Q$-{\it CY}
variety
if $W$ is a rational curve, an abelian variety,
a $\Q$-torus, or a $\Q$-CY variety, respectively.

Here a projective variety is called a {\it $\Q$-torus} if
it has an abelian variety as its finite \'etale (Galois)
cover. A projective variety $X$ with only klt singularities
is called a $\Q$-{\it CY variety} if some positive multiple of
its canonical divisor $K_X$ is linearly equivalent to
the trivial divisor. A simply connected
$\Q$-{\it CY} surface with only Du Val
singularities is called a {\it normal K3 surface}.
A $\Q$-torus is an example of a smooth $\Q$-CY.

\begin{theorem} (cf.~\cite[Theorem 1.4]{LZ}) \label{ThF'}
Let $X$ be a smooth projective variety of dimension $n$
and $D$ a reduced
divisor on $X$ with simple normal crossings.
Then the closure $\NE(X)$ of effective $1$-cycles on $X$ is generated by
the $(K_X + D)$-non negative part and
at most a countable collection of
extremal $\mathcal{D}$-rational curves
$\{\ell_i\}_{i \in N}, \, N \subseteq \Z$:
$$\NE(X) = \NE(X)_{K_X + D \ge 0} +
\sum_{i \in N} \R_{> 0}[\ell_i] .$$
Further, $-\ell_i . (K_X + D)$ is in $\{1, 2, \dots, 2n\}$,
and is equal to $1$
when $\ell_i$ is not a $\mathcal{D}$-{\rm compact} rational curve.
In particular, if the pair $(X, D)$ is Mori hyperbolic
then $K_X + D$ is nef.
\end{theorem}

\begin{theorem} (cf.~\cite[Theorem 1.5]{LZ}) \label{ThE'}
Let $X$ be a smooth projective variety of dimension $n$
and $D$ a reduced
divisor on $X$ with simple normal crossings.
Fix an $r$ in $\{1, 2, \dots, n\}$.
Assume that {\rm Abundance($l$)} holds
for $l \leq r$  and that
$D$ has at least $n-r+1$
irreducible components amongst which at least $n-r$ are ample.
If $K_X+D$ is not ample, then either
it has non-positive degree on a $\mathcal{D}$-rational curve
or on a $\mathcal{D}$-algebraic $1$-torus,
or some positive multiple of
it is linearly equivalent to the trivial divisor on
a smooth $\mathcal{D}$-$\Q$-CY
variety $T$ of dimension $< r$.
We can take $T$ to be a $\Q$-torus if
further {\rm CY($m$)} holds for $m<r$.

In particular, if {\rm Abundance($l$)}
and {\rm CY($m$)} hold for $l \leq r$ and $m<r$
and the pair $(X,D)$ is Brody hyperbolic then $K_X+D$ is ample.
\end{theorem}

Since the two conjectures are known for $r = 3$ (see the
remark above),
an immediate corollary is the following generalization,
to arbitrary dimensions for pairs, of results that were only known
up to dimension two in the case $D=0$ concerning hyperbolic
varieties.

\begin{theorem} (cf.~\cite[Theorem 1.6]{LZ}) \label{ThEE''}
Let $X$ be a smooth projective variety of dimension $n$
and $D$ a reduced
divisor on $X$ with simple normal crossings such that
$D$ has at least $n-2$
irreducible components amongst which at least $n-3$ are ample.
If $K_X+D$ is not ample, then either
it has non-positive degree on a $\mathcal{D}$-rational curve
or on a $\mathcal{D}$-algebraic $1$-torus,
or some positive multiple of
it is linearly equivalent to
the trivial divisor on
a $\mathcal{D}$-torus
$T$ with $\dim T \le 2$.
In particular, if $(X,D)$ is Brody hyperbolic, then $K_X+D$ is ample.
\end{theorem}

We remark that Theorems \ref{ThF'}, \ref{ThEE''} and \ref{ThE'}
are special cases of Theorems \ref{ThF} and \ref{ThE} (see also Theorem \ref{ThD})
where we allow the pair
$(X,D)$ to be singular%
\footnote{Here by a singular pair, we will assume that it is
a dlt pair.
The assumption is natural (and in many respects
the most general)
as our proof is by induction on dimension
from running the LMMP for singular pairs. It implies
an explicit adjunction formula
and
that $k$-fold intersections of components of $D$
are of pure codimension-$k$ in $X$, crucial
in our inductive procedure.}%
$\ $
and $D$ to be
augmented with a fractional
divisor
$\Gamma$ in a setting that is natural to our approach.
Theorem \ref{ThF'} is also obtained in
McQuillan-Pacienza \cite{MP} by a different method, but not our more
general result in the singular case, Theorem \ref{ThF}.
That article obtained the results in the setting of
complete intersection ``stacks'' by a direct analysis
of Mori's bend and break procedure and includes
for the most part our refined geometric version of the cone theorem
for smooth pairs, Theorem~\ref{ThF'}, which we obtain
however in our
more general singular (dlt) setting, Theorem~\ref{ThF}, and
by a different and independent method.

By running the minimal model program, an
easy consequence of this cone theorem is the following criterion for
$K_X+D$ to be peudo-effective (the weakest form of
positivity in birational geometry) under a mild condition on $D$.
It is a special case of Theorem \ref{ThH}.

\begin{theorem} (cf.~\cite[Theorem 1.7]{LZ}) \label{ThH'}
Let $X$ be a smooth projective variety
and $D$ a nonzero simple normal crossing reduced divisor
no component of which is uniruled.
Then $K_X + D$ is not pseudo-effective if and only if
$D$ has exactly one irreducible component and
$X$ is dominated by $\mathcal D$-rational curves $\ell$
with
$$-\ell . (K_X + D) = 1, \,\,\,\, \ell \cong \PP^1, \,\,\,\,
\ell \setminus \ell \cap D \cong \C$$
for a general $\ell$.

In particular, $K_X+D$ is pseudo-effective for a log
smooth pair $(X,D)$ such that
$D$ contains two or more non-uniruled components.
\end{theorem}

Key to our proof in the presence of Cartier boundary divisors
is Kawamata's result
on the length of extremal rays, a fundamental
result in the subject.
Our proofs of Theorems $\ref{ThE} \sim \ref{ThD}$ are inductive in nature and
reduce the problem to questions on
adjoint divisors
in lower dimensions by adjunction.
The log minimal model program (LMMP) is run formally without going into
its technical details.
Our inductive procedure is naturally adapted to answer some
fundamental questions concerning adjoint divisors.

Theorems \ref{ThF}, \ref{ThE} and \ref{ThH} below
include Theorems \ref{ThF'}, \ref{ThE'} and \ref{ThH'} as special cases.\\[-3mm]

Let $X$ be a normal projective variety and $D$ a reduced Cartier divisor on $X$.
The pair $(X, D)$ is called {\it BH, or MH with respect to a $($Cartier$)$
decomposition} of $D = \sum_{i=1}^s D_i$,
if  each $D_i$ is a reduced
(Cartier) divisor on $X$, not necessarily irreducible,
and both $X \setminus D$ and $(\cap_{i \in I} D_i) \setminus (\cup_{j \in J} D_j)$
are respectively BH, or MH
for every partition of $\{1, \dots, s\} = I \coprod J$.
Note that, for a log smooth pair $(X, D)$,
the respective definitions of hyperbolicity given early on
are equivalent to the
definitions of hyperbolicity here with respect to the irreducible
decomposition of $D$.

\begin{theorem} (cf.~\cite[Theorem 3.1]{LZ}) \label{ThF}
Let $X$ be a
projective variety of dimension $n$,
$D$ a reduced divisor that is the sum of a
collection $\mathcal D=\{D_i\}_{i=1}^s$ of reduced
Cartier divisors
and $\Gamma$
an effective Weil $\Q$-divisor on $X$ such that the pair $(X, D + \Gamma)$
is divisorial log terminal.
Then the closure $\NE(X)$ of effective $1$-cycles on $X$ is generated by
the $(K_X + D + \Gamma)$-non negative part and
at most a countable collection of
extremal $\mathcal{D}$-rational curves
$\{\ell_i\}_{i \in N}, \, N \subseteq \Z$:
$$\NE(X) = \NE(X)_{K_X + D + \Gamma \ge 0} +
\sum_{i \in N} \R_{> 0}[\ell_i] .$$
Furthermore,
$-\ell_i . (K_X + D+\Gamma)$ is in $(0, 2n]$,
and is even in $(0, 2)$
if $\ell_i$ is not $\mathcal{D}$-{\rm compact}.
In~particular, if the pair $(X,D)$ is Mori hyperbolic with respect to the
Cartier decomposition $D = \sum_{i=1}^s D_i$ then $K_X+D+\Gamma$ is nef.
\end{theorem}

One may take $\Gamma = 0$ in Theorems $\ref{ThF} \sim \ref{ThD}$.
Such a $\Gamma$ naturally appears in our inductive procedure as the `different'
in the adjunction formula.

We remark that {\rm Abundance($l$)} (for dlt pairs)
and {\rm CY($m$)}
always hold for $l \leq 3$ and $m\leq 2$ (\cite[\S 3.13]{KM}, \cite{MM})
and that a normal K3 surface has infinitely many
elliptic curves (\cite{MM}).

\begin{theorem} (cf.~\cite[Theorem 3.2]{LZ}) \label{ThE}
Let $X$ be a
projective variety of dimension $n$,
$D$ a reduced Weil divisor and $\Gamma$
an effective Weil $\Q$-divisor on $X$ such that the pair $(X, D + \Gamma)$
is divisorial log terminal $($dlt$)$.
Assume one of the following three conditions.
\begin{itemize}
\item[(1)]
$n \le 2$.
\item[(2)]
$n = 3$, $D$ is nonzero and Cartier.
\item[(3)]
$n \ge 4$, the pair $(X, D + \lfloor \Gamma \rfloor)$ is log smooth,
there is an $r \in \{1, 2, \dots, n\}$ such that
$D$ has at least $n-r+1$
irreducible components amongst which at least $n-r$ are ample
and {\rm Abundance($l$)} $($for dlt pairs$)$
holds for $l \leq r$.
\end{itemize}
Then $K_X+D+\Gamma$ is ample, unless either
$(K_X + D + \Gamma)$ has non-positive degree on a
$\mathcal{D}$-rational curve or on a $\mathcal{D}$-algebraic $1$-torus,
or the restriction $(K_X + D + \Gamma)_{| T} \sim_{\Q} 0$
for a $\mathcal{D}$-$\Q$-CY variety $T$.
When $T$ is singular, it is
a normal $K3$ surface and $n\in \{2, 3\}$.
%
Otherwise, $T$ can be taken to be an abelian surface
or an elliptic curve in Cases~$(1)$,~$(2)$ and,
if further {\rm CY($m$)} holds for $m<r$,
taken to be a $\Q$-torus with $ \dim T < r$
in Case~$(3)$.

In particular, if the pair $(X, D)$ is Brody hyperbolic
then $K_X+D+\Gamma$ is ample, provided that either $n \le 3$ or
{\rm CY($m$)} holds for $m<r$.
\end{theorem}

\begin{theorem} (cf.~\cite[Theorem 3.3]{LZ}) \label{ThD}
Let $X$ be a
projective variety of dimension $n$,
$D$ a reduced divisor that is the sum of a
collection $\mathcal D=\{D_i\}_{i=1}^s$ of reduced
Cartier divisors for some $s \ge 1$ and
$\Gamma$
an effective Weil $\Q$-divisor on $X$ such that the pair $(X, D + \Gamma)$
is divisorial log terminal.
Assume one of the following two conditions.
\begin{itemize}
\item[(1)]
$n = 4$, $s \ge 2$ and $D_1$ is irreducible and ample.
\item[(2)]
$n \ge 4$, $s \ge n-r +1$ for some $r \in \{1, 2, \dots, n\}$,
$D_j$ is ample for all $j \le n-r+1$ and {\rm Abundance($l$)}
$($for dlt pairs$)$ holds for $l \leq  r$.
\end{itemize}
Then $K_X + D + \Gamma$ is ample, unless
either $(K_X + D + \Gamma)$ has non-positive degree on a
$\mathcal{D}$-rational curve or on a $\mathcal{D}$-algebraic $1$-torus, or
the restriction
$(K_X + D + \Gamma)_{| T} \sim_{\Q} 0$
for some $\mathcal{D}$-$\Q$-CY
variety $T$ $( \dim T < r$ in Case$(2))$.
$T$ is smooth when
the pair
$(X, D + \lfloor \Gamma \rfloor)$ is log smooth
and can be taken to be a $\Q$-torus when further {\rm CY($m$)} holds for $m<r$.
\end{theorem}

%
%

\begin{theorem} (cf.~\cite[Theorem 3.4]{LZ}) \label{ThH}
Let $X$ be a $\Q$-factorial projective variety,
and $D$ a nonzero reduced Weil divisor with $D = \sum_{i=1}^s D_i$ the irreducible decomposition
such that the pair $(X,D)$ is divisorial log terminal.
\begin{itemize}
\item[(1)]
Assume that no component of $D$ is uniruled
and $K_X + D$ is not pseudo-effective.
Then $D$ has exactly one irreducible component
and $X$ is dominated by $\mathcal D$-rational curves $\ell$
with
$$-\ell . (K_X + D) = 1, \,\,\,\, \ell \cong \PP^1, \,\,\,\, \ell \setminus \ell \cap D \cong \C$$
for a general $\ell$.
\item[(2)]
Conversely, assume that $X$ has at worst canonical singularities and
$X$ is dominated by $\mathcal D$-rational curves. Then $K_X + D$ is not pseudo-effective.
\end{itemize}
\end{theorem}

\noindent
Of interest here is that the above
geometric criterion for the pseudo-effectivity of $K_X+D$
in Theorem~\ref{ThH}
is obtained naturally from Theorem~\ref{ThF}
by running the log minimal model program for the pair $(X,D)$.
In fact, by running the log minimal model program for a normal surface
we get the following slightly stronger consequence.
It is noteworthy that Keel-${\text M}^{\text c}$Kernan \cite{KMc}
has obtained the same
conclusion in dimension two
without the assumption of the absence of $\mathcal D$-rational curves in $D$.

\begin{proposition} (cf.~\cite[Proposition 3.5]{LZ}) \label{PropE}
Let $X$ be a projective surface and $D$ a nonzero reduced
Weil divisor with $D = \sum_{i=1}^s D_i$ the irreducible decomposition
such that the pair $(X, D)$ is divisorial log terminal
and there are no
$\mathcal D$-rational curves in $D$.
Assume that $K_X + D$ is not pseudo-effective.
Then $D$ has exactly one irreducible component,
and $X$ is dominated by $\mathcal D$-rational curves $\ell$
with
$$-\ell . (K_X + D) = 1, \,\,\,\, \ell \cong \PP^1, \,\,\,\, \ell \setminus \ell \cap D \cong \C$$
for a general $\ell$.
\end{proposition}

\noindent
{\bf N.B.} In Theorem \ref{ThD} (1), $D_1$ needs to be irreducible
so that Theorem \ref{ThE} (2)  can be applied:
if $D_1 = D_{11} + D_{12}$ is reducible
then the restriction $D_{2 |D_{11}}$ might be zero.
{\it The ampleness assumption on
some $D_i$'s in various theorems is due to the same reason.}

\section{Positivity of canonical divisors of singular varieties and Lang hyperbolicity}

This section is based on the paper \cite{SingHyp}.
Let $X$ be a variety.
$X$ is {\it Brody hyperbolic} (resp. {\it algebraic Lang hyperbolic})
if every holomorphic map $V \to X$ where $V$ is the complex line $\C$ (resp. $V$ is
an abelian variety) is a constant map. Since an abelian variety is a complex torus,
Brody hyperbolicity implies algebraic Lang hyperbolicity.
When $X$ is compact complex variety, Brody hyperbolicity is equivalent to the
usual Kobayashi hyperbolicity (cf. \cite{La}).

In the first part (Thoerem \ref{ThA5})
of this section,
we let $X$ be a normal projective variety and
aim to show the ampleness of the
{\it canonical divisor $K_X$ of $X$}, assuming that
$X$ is algebraic Lang hyperbolic. We allow $X$ to have arbitrary singularities and
assume only that $X$ is {\it $\Q$-Gorenstein} (so that the ampleness of $K_X$ is well-defined),
i.e., $K_X$ is {\it $\Q$-Cartier:}
$mK_X$ is a Cariter divisor for some positive integer $m$.

For related work,
it was proved in \cite{Pe} that a $3$-dimensional
hyperbolic smooth projective variety $X$ has ample $K_X$ unless $X$ is
a Calabi-Yau manifold where every non-zero effective divisor is ample.
The authors of \cite{HLW} prove the ampleness of $K_X$ when $X$ is a smooth
projective threefold having a K\"ahler metric of negative holomorphic sectional curvature;
they also generalize the results to higher dimensions with some additional conditions.

In the second part of this section (Theorems \ref{ThC5})
aiming to prove Lang's conjecture in Corollary \ref{Cor3},
even the normality of $X$ is not assumed.
Our approach is to take a projective resolution of $X$ and run the relative Minimal Model Program
(MMP) over $X$.
We use only the frame work of MMP, but not its detailed technical part.
Certain mild singularities occur naturally along the way.
See \cite[Definition  2.34]{KM} for definitions of {\it canonical},
{\it Kawamata log terminal} (klt), and {\it divisorial log terminal} (dlt) {\it singularities}.

In the last part
(Proposition \ref{PropA}), we try to avoid assuming conjectures.

Now let $X$ be a normal projective variety.
$X$ is a {\it Calabi-Yau variety} if
$X$ has at worst canonical singularities,
its canonical divisor is $\Q$-linearly equivalent to zero: $K_X \sim_{\Q} 0$
(and hence $X$ has Kodaira dimension $\kappa(X) = 0$)
and the {\it irregularity} $q(X) := h^1(X, \OO_X) = 0$.
$X$ is {\it of general type} if some (equivalently every) projective resolution $X'$ of
$X$ has maximal Kodaira dimension: $\kappa(X') = \dim X'$.
A $\Q$-Gorenstein variety $X$ is {\it minimal} if the canonical divisor $K_X$ is {\it numerically
effective} (= {\it nef}).

Conjecture \ref{nALH} below is long standing.
When $\dim X \le 2$, it is true by the classification of
complex surfaces and the {\it fact $(*)$}: a (smooth) $K3$ surface has infinitely many
(singular) elliptic curves (cf. \cite[Theorem in Appendix]{MM}).

In Conjecture \ref{nALH}, the conclusion means the existence of at least one non-constant
holomorphic map $f : V \to X$ from an abelian variety $V$,
but does {\it not} require the image $f(V)$ (or the union of such images) to be Zariski-dense
in $X$. This does not seem sufficient for our purpose to show the non-existence of subvariety $X'$ of Kodaira dimension
zero in an algebraic Lang hyperbolic variety $W$ as in Corollary \ref{Cor3} below
(see \ref{setup1.1}, and think about a proof of the non-hyperbolicity of every normal $K3$ surface
using the fact $(*)$ above).
Fortunately, we are able to show in Corollary \ref{Cor3}
(or Theorems \ref{ThC5})
that the normalization $X$ of $X' \subseteq W$ is a Calabi-Yau variety
and hence $f$ composed with the finite morphism $X \to X' \subseteq W$
produces a non-constant holomorphic map from the abelian variety $V$,
thus deducing a contradiction to the hyperbolicity of $W$.

\begin{conjecture}\label{nALH}
Suppose that $X$ is an absolutely minimal {\rm (in the sense of \ref{rThA5})} Calabi-Yau variety.
Suppose further that every birational morphism $X \to Y$ onto a normal projective variety is an isomorphism.
Then $X$ is not algebraic Lang hyperbolic.
\end{conjecture}

We recall the definition of the {\it nef reduction} $f_L : X \dasharrow Y$
of a nef $\Q$-Cartier divisor $L$ on a normal projective variety
(cf. \cite{8aut}).
$f_L$ is almost holomorphic with $F$ a compact general fibre;
the restriction $L_{| F}$
is numerically trivial;
a general curve $C$ on $X$
is contracted by $f_L$ if and only if $L . C = 0$; this $f_L$ is unique up to
birational equivalence of $Y$; $\dim Y$ is called the {\it nef dimension} of $L$
and denoted as $n(L)$. We have $n(L) = \dim X$ (maximal case) if and only if
we have the positivity: $L . C > 0$ for every curve $C$ which is not contained in a fixed countable union
of Zariski-closed proper subsets of $X$; also $n(L) = 0$ holds
if and only if $L \equiv 0$ (numerically zero).

The abundance conjecture below is slightly weaker than the usual one where
the extra hypothesis $(*)$ is usually not assumed.
When $K_X$ is nef and big or when
$\dim X \le 3$, both Conjecture \ref{Abund}(1) $\sim$ \ref{Abund}(2) (and their log version, even
without the hypothesis $(*)$) are true; see \cite[Theorem 3.3, \S 3.13]{KM}.

\begin{conjecture}\label{Abund}
Let $X$ be an $n$-dimensional minimal normal projective variety, i.e.,
the canonical divisor $K_X$ is a nef $\Q$-Cartier divisor.
Assume the hypothesis $(*)$: the nef dimension $n(K_X)$ {\rm (cf.~\cite{8aut})} satisfies $n(K_X) = n$.
\begin{itemize}
\item[(1)]
If $X$ has at worst klt singularities,
then
$K_X$ is semi-ample, i.e.,
the linear system $|mK_X|$ is base-point free for some $m > 0$.
\item[(2)]
If $X$ has at worst canonical singularities and
$K_X \not\equiv 0$ $($not numerically zero$)$, then
the Kodaira dimension $\kappa(X) > 0$.
\end{itemize}
\end{conjecture}

{\it In this paper, by hyperbolic we mean algebraic Lang hyperbolic.}

\par \vskip 1pc
When $X$ has at worst klt singularities, Theorem \ref{ThA5} below follows from
the MMP and has been generalized to the quasi-projective case
in \cite{LZ}.

In Theorem \ref{ThA5}, we don't impose any condition on the
singularities of $X$, except the $\Q$-Cartierness of $K_X$ for the concluding
amplenss of $K_X$ to make sense.

Without assuming Conjecture \ref{nALH} or \ref{Abund}, we can
at least say that $K_{X}$ is movable or nef in codimension-one
(cf. Remark \ref{rThA5}).

\begin{theorem} (cf.~\cite{SingHyp}) \label{ThA5}
Let $X$ be a $\Q$-Gorenstein normal projective variety which is algebraic Lang hyperbolic.
Assume that Conjectures $\ref{nALH}$ and $\ref{Abund}(1)$ hold. Then
$K_X$ is ample at smooth points and klt points of $X$.
To be precise, there is a birational morphism $f_c : X_c \to X$
such that $X_c$ has at worst klt singularities, $K_{X_c}$ is ample, and
$$E_c := f_c^*K_X - K_{X_c}$$
is effective and $f_c$-exceptional
with $f_c(E_c) \subseteq \Nklt(X)$, the non-klt locus of $X$
{\rm (cf. \cite[\S 4.4]{Fj}).} In particular, $f_c = \id$ and $K_X$
is ample, if $X$ has at worst klt singularities {\rm (cf. \ref{rThA5})}.
\end{theorem}

Even the normality of $X$ is not assumed here in Theorem \ref{ThC5}.

\begin{theorem} (cf.~\cite{SingHyp}) \label{ThC5}
Let $W$ be an algebraic Lang hyperbolic projective variety, and
$X \subseteq W$ a projective subvariety $($with $X = W$ allowed$)$.
Assume either $\dim X \le 3$ or Conjecture $\ref{Abund}(2)$
$($resp. either $\dim X \le 3$ or Conjecture $\ref{Abund}(2)$ with the hypothesis $(*)$ removed$)$.
Then there is a birational surjective morphism $g_m : X_m \to X$ such that
$X_m$ is a minimal variety with at worst canonical
singularities and one of the following is true.
\begin{itemize}
\item[(1)]
$K_{X_m}$ is ample. Hence both $X_m$ and $X$ are of general type.
\item[(2)]
$g_m : X_m \to X$ is the normalization map; $X_m$ is an absolutely minimal
{\rm (in the sense of \ref{rThA5})} Calabi-Yau variety with $\dim X_m \ge 3$.
\item[(3)]
There is an almost holomorphic map $\tau : X_m \dasharrow Y$
$($resp. a holomorphic map $\tau: X_m \to Y)$
such that its general fibre $F$ is an absolutely minimal
Calabi-Yau variety with $3 \le \dim F < \dim X_m$,
and $(g_m)_{|F} : F \to g_m(F) \subset X$ is the normalization map.
\end{itemize}
\end{theorem}

Given a projective variety $W$, let $\widetilde{W} \overset{\sigma}\to W$ be a
projective resolution. We define the {\it albanese variety
of $W$} as $\Alb(W) := \Alb(\widetilde{W})$, which is independent of the choice
of the resolution $\widetilde{W}$, since every two resolutions of $W$ are
dominated by a third one and the albanese variety, being an abelian variety,
contains no rational curves. We define
the {\it albanese $($rational$)$ map} $\alb_W : W \dasharrow \Alb(W)$
as the composition
$$W \overset{\sigma^{-1}}\dasharrow \widetilde{W}
\overset{\alb_{\widetilde{W}}}\longrightarrow \Alb(\widetilde{W}) .$$

Lang's \cite[Conjecture 5.6]{La} (one direction only; see \ref{rThA5} for the other direction)
follows from Conjectures \ref{nALH} and \ref{Abund}(2):

\begin{corollary} (cf.~\cite{SingHyp}) \label{Cor3}
Let $W$ be an algebraic Lang hyperbolic projective variety.
Assume either $\dim W \le 3$ or Conjecture $\ref{Abund}(2)$. Then we have:
\begin{itemize}
\item[(1)]
If Conjecture $\ref{nALH}$ holds, then $W$ and all its subvarieties are of general type.
\item[(2)]
If the albanese map  $\alb_W : W \dasharrow \Alb(W)$
has general fibres
of dimension $\le 2$, then $W$ is of general type.
\end{itemize}
\end{corollary}

If we don't assume Conjecture \ref{nALH} in Theorem \ref{ThA5},
then the proof shows:

\begin{corollary} (cf.~\cite{SingHyp})
Let $X$ be a $\Q$-Gorenstein normal projective variety of dimension $n$, which is algebraic Lang hyperbolic.
Assume either $n \le 3$ or Conjecture $\ref{Abund}(1)$. Then
either $K_X$ is ample at smooth points and klt points of $X$ as
detailed in Theorem $\ref{ThA5}$;
or $X$ is a Calabi-Yau variety of dimension $n \ge 3$; or
$X$ is covered by subvarieties $\{F_t'\}$ whose normalizations
are Calabi-Yau varieties of the same dimension $k$ with $3 \le k < n$.
\end{corollary}

Without assuming Conjectures \ref{Abund} (and \ref{nALH}), we have:

\begin{proposition} (cf.~\cite{SingHyp}) \label{PropA}
Let $W$ be a projective variety
which is algebraic Lang hyperbolic. Then $W$ is of general type
provided that either the albanese map $\alb_W : W \dasharrow \Alb(W)$
is generically finite $($but not necessarily surjective$)$; or
the Kodaira dimension $\kappa(W) \ge \dim W - 2$;
or $\kappa(W) \ge \dim W - 3$ and Conjecture $\ref{nALH}$ holds in dimension three.
\end{proposition}

\begin{remark}\label{rThA5}
(1) In Theorem \ref{ThA5},
by the equality $f_c^*K_X = K_{X_c} + E_c$ and the ampleness of $K_{X_c}$, the {\it exceptional locus} $\Exc(f_c)$ (the subset of $X_c$ along which
$f_c$ is not isomorphic) is contained in $\Supp E_c$. This and the effectivity of $E_c$ justifies the
phrasal: $K_{X_c}$ is ample outside $f_c(E_c)$.

(2) Without assuming Conjecture \ref{nALH} or \ref{Abund},
the proof of Theorem \ref{ThA5}
shows that $f^*K_X = K_{X'} + E'$ with $K_{X'}$ nef and $E'$ $f$-exceptional.
Hence $K_X = f_*K_{X'}$ is movable, or nef in codimension-one.

(3) Let $X_2 \to X_1$ be a finite morphism (but not necessarily surjective).
If $X_1$ is Brody hyperbolic or algebraic Lang hyperbolic then so is $X_2$.
The converse is not true.

(4) Every algebraic Lang hyperbolic projective variety $X_1$ is {\it absolutely minimal} in the sense:
if $X_2$ is a normal projective variety with at worst klt singularities then every birational map
$h: X_2 \dasharrow X_1$ is a well defined morphism. (This result was proved by S.~Kobayashi when $X_2$
is nonsingular). Indeed, let $Z$ be a resolution of the graph of $h$ such that
we have birational surjective morphisms $p_i : Z \to X_i$
satisfying $h \circ p_2 = p_1$. Then every fibre $p_2^{-1}(x_2)$ is rationally chain connected
by \cite[Corollary 1.5]{HM} and hence $p_1(p_2^{-1}(x_2))$ is a point since
hyperbolic $X_1$ contains no rational curve.
Thus $h$ is extended to a well defined morphism by \cite[Proof of Lemma 14]{Ka81},
noting that $X_2$ is normal and $p_2$ is surjective, and using the Stein factorization.

(5) If $Y$ is an algebraic Lang hyperbolic Calabi-Yau variety (like $X_m$ and $F$ in
\ref{ThC5} (2) and (3), respectively), then
every birational morphism $h: Y \to Z$ onto a normal projective variety is an isomorphism.
Indeed, by \cite[Corollary 1.5]{Ka88}, $Z$ has only
canonical singularities. Thus the exceptional locus $\Exc(h)$ is covered by
rational curves by \cite[Corollary 1.5]{HM}. Since $Y$ is hyperbolic and hence has no rational
curve, we have $\Exc(h) = \emptyset$ and hence $h$ is an isomorphism, $Z$ being normal
and by the Stein factorization.

(6) About the converse of Corollary \ref{Cor3}, i.e., another direction of Lang
\cite[Conjecture 5.6]{La}:
 Assume that every non-uniruled projective variety has a minimal model with at worst
canonical singularities
and that abundance Conjecture \ref{Abund}(2) holds.
If a projective variety $W$ and all its subvarieties are of
general type, then $W$ is algebraic Lang hyperbolic. Indeed, let $f: V \to W$
be a morphism from an abelian variety $V$
and let $V \to X \to f(V)$ be its Stein factorization,
where $V \to X$ has connected general fibre $F$ and $X \to f(V)$
is a finite morphism. Since $V$ is non-uniruled, so is $F$. Hence $\kappa(F) \ge 0$ by the assumption.
The assumption and Iitaka's $C_{n,m}$ also
imply that
$$0 = \kappa(V) \ge \kappa(F) + \kappa(X) \ge \kappa(X) \ge \kappa(f(V)) = \dim f(V)$$
(cf.~\cite[Corollary 1.2]{Ka85}). Hence $f$ is a constant map.
\end{remark}

\begin{setup}\label{setup1.1}
{\it Comments about the proofs:}
{\rm
In our proofs, neither the existence of minimal model nor the termination
of MMP is assumed.
To show that every subvariety $X$ of a hyperbolic $W$ is of general type,
one key observation is that a relative minimal model of $X$
birational over $X$ which exists by the main Theorem $1.2$ in \cite{BCHM},
{\it is indeed a minimal model of $X$} since $W$ is hyperbolic.
One natural approach is to take a general fibre $F$ (which may not even be normal)
of an Iitaka (rational)
fibration of $X$ (assuming $\kappa(X) \ge 0$)
and prove that $F$ has a minimal model $F_m$. Next, one tries to show that $q(F_m) = 0$
and $F_m$ is a Calabi-Yau variety and then tries to use Conjecture \ref{nALH} to produce
a non-hyperbolic subvariety $S$ of $F_m$, but this does not guarantee the
same on $F \subset X$ (to contradict the hyperbolicity of $X$)
because such $S \subseteq F_m$ might be contracted to a point on $F$.
In our approach, we are able to show that the {\it normalization of $F$ is a Calabi-Yau variety,
which is the key of the proofs.} It wouldn't help even one assumes the smoothness of the ambient space
$W$ since its subvariety $X$ may not be smooth, or at least normal or Cohen-Macaulay to define
the canonical divisor $K_X$ meaningfully to pull back or push forward.
}
\end{setup}

\end{document}